\documentclass[12pt]{article}
\usepackage{amsmath}
\usepackage{amssymb}
\usepackage{amsthm}
\usepackage[mathscr]{eucal}
\usepackage{latexsym}
\usepackage{amscd}

\begin{document}

\voffset -1.5truecm
\oddsidemargin0.5truecm
\evensidemargin0.5truecm

\theoremstyle{plain}
\newtheorem{lema}{Lemma}[section]
\newtheorem{prop}[lema]{Proposition}
\newtheorem{coro}[lema]{Corolary}
\newtheorem{teor}[lema]{Theorem}
\newtheorem{demo}[lema]{Proof of Theorem 1}
\newtheorem{ejem}[lema]{Example}
\newtheorem{defi}[lema]{Definition}
\newtheorem{obse}[lema]{Remark}
\newtheorem{obses}[lema]{Remarks}

\renewcommand{\refname}{\large\bf References}
\renewcommand{\thefootnote}{\fnsymbol{footnote}}

\def\flecha{\longrightarrow}
\def\asocia{\longmapsto}
\def\sii{\Longleftrightarrow}
\def\vector#1#2{({#1}_1,\dots,{#1}_{#2})}
\def\conjunto#1{\boldsymbol{[}\,{#1}\,\boldsymbol{]}}
\def\particion{\vdash}

\def\natural{{\mathbb N}}
\def\entero{{\mathbb Z}}
\def\racional{{\mathbb Q}}
\def\real{{\mathbb R}}
\def\complejo{{\mathbb C}}

\def\cara#1{\chi^{#1}}
\def\permu#1{\phi^{\,#1}}
\def\coefi{{\sf k}(\la,\mu,\nu)}
\def\coefili#1#2#3{{\sf k}({#1},{#2},{#3})}
\def\kron{\chi^\lambda\otimes\chi^\mu}
\def\kronli#1#2{\chi^{#1}\otimes\chi^{#2}}

\def\gruposim#1{{\sf S}_{#1}}

\def\la{\lambda}
\def\longi#1{\ell({#1})}

\def\menore{\prec}
\def\menori{\preccurlyeq}
\def\mayore{\succ}
\def\mayori{\succcurlyeq}

\def\abc#1#2#3{{#1}_{{#2}{#3}}}
\def\columna#1{{\sf col}({#1})}
\def\contenido#1{{\sf cont}({#1})}
\def\forma#1{{\sf sh}({#1})}
\def\cano#1{{\sf C}({#1})}
\def\insertar#1#2#3{{#1}_{#2} \rightarrow ( \cdots ( {#1}_2 \rightarrow
({#1}_1 \rightarrow {#3})) \cdots )}
\def\kostka#1#2{K_{{#1}{#2}}}
\def\kostkasf#1#2{{\sf K}_{{#1}{#2}}}

\def\lr{{\sf LR}(\alpha,\beta;\nu)}
\def\lrest{{\sf LR}^*(\alpha,\beta;\nu)}
\def\lrd{{\sf LR}^*(\alpha,\beta;\nu)}

\def\matriz{{\sf M}(\lambda,\mu)}
\def\numatriz{{\sf m}(\lambda,\mu)}
\def\matrizcon{{\sf M}_{\nu}(\lambda,\mu)}
\def\cardimatrizcon{{\sf m}_\nu(\lambda,\mu)}
\def\partipla{{\sf P}_\nu(\la,\mu)}
\def\matrizd{{\sf M}^*(\lambda,\mu)}
\def\numatrizd{{\sf m}^*(\lambda,\mu)}
\def\matriztres{{\sf M}(\lambda,\mu,\nu)}
\def\numatres{{\sf m}(\lambda,\mu,\nu)}
\def\matriztresest{{\sf M}^*(\lambda,\mu,\nu)}
\def\matriztresd{{\sf M}^*(\lambda,\mu,\nu)}
\def\numatresd{{\sf m}^*(\lambda,\mu,\nu)}
\def\matriztresdc{{\sf M}^*(\lambda,\mu,\nu^{\,\prime})}
\def\parti#1{\pi(#1)}
\def\liri#1#2#3{{\sf lr}({#1},{#2};{#3})}
\def\lirid#1#2#3{{\sf lr}^*({#1},{#2};{#3})}

\def\palabra#1#2#3{{#1}_{#2}\cdots{#1}_{#3}}
\def\renglon#1{{\sf row}(#1)}
\def\sst{semistandard\ tableau\ }
\def\sucesion#1#2{{#1}_1 < {#1}_2 < \cdots < {#1}_{#2}}
\def\wcol#1{w_{\rm col}({#1})}

\def\bili#1#2{\langle{#1}, {#2}\rangle}


\setlength\unitlength{0.08em}
\savebox0{\rule[-2\unitlength]{0pt}{10\unitlength}%
\begin{picture}(10,10)(0,2)
\put(0,0){\line(0,1){10}}
\put(0,10){\line(1,0){10}}
\put(10,0){\line(0,1){10}}
\put(0,0){\line(1,0){10}}
\end{picture}}

\newlength\cellsize \setlength\cellsize{18\unitlength}
\savebox2{%
\begin{picture}(18,18)
\put(0,0){\line(1,0){18}}
\put(0,0){\line(0,1){18}}
\put(18,0){\line(0,1){18}}
\put(0,18){\line(1,0){18}}
\end{picture}}
\newcommand\cellify[1]{\def\thearg{#1}\def\nothing{}%
\ifx\thearg\nothing
\vrule width0pt height\cellsize depth0pt\else
\hbox to 0pt{\usebox2\hss}\fi%
\vbox to 18\unitlength{
\vss
\hbox to 18\unitlength{\hss$#1$\hss}
\vss}}
\newcommand\tableau[1]{\vtop{\let\\=\cr
\setlength\baselineskip{-16000pt}
\setlength\lineskiplimit{16000pt}
\setlength\lineskip{0pt}
\halign{&\cellify{##}\cr#1\crcr}}}
\savebox3{%
\begin{picture}(15,15)
\put(0,0){\line(1,0){15}}
\put(0,0){\line(0,1){15}}
\put(15,0){\line(0,1){15}}
\put(0,15){\line(1,0){15}}
\end{picture}}
\newcommand\expath[1]{%
\hbox to 0pt{\usebox3\hss}%
\vbox to 15\unitlength{
\vss
\hbox to 15\unitlength{\hss$#1$\hss}
\vss}}


\begin{centering}
{\Large\bf  Kronecker products and the RSK correspondence}\\[.8cm]
{{\large\sf Diana Avella-Alaminos \footnotemark[1]}}\\[.1cm]
Universidad Nacional Aut\'onoma de M\'exico\\
Facultad de Ciencias, Departamento de Matem\'aticas\\
04510 M\'exico, D.F., MEXICO\\
{\tt avella@matem.unam.mx}\\[.3cm]
{\large\sf and}\\[.3cm]
{{\large\sf Ernesto Vallejo \footnotemark[2]}\footnotemark[3]}\\[.1cm]
Universidad Nacional Aut\'onoma de M\'exico\\
Instituto de Matem\'aticas, Unidad Morelia\\
Apartado Postal 61-3, Xangari\\
58089 Morelia, Mich., MEXICO\\
{\tt vallejo@matmor.unam.mx}\\
\end{centering}
\footnotetext[1]{Supported by UNAM-DGAPA, UNAM-DGEP and CONACYT Mexico.}
\footnotetext[2]{Corresponding author.}
\footnotetext[3]{Supported by CONACYT Mexico grant 47086-F and UNAM-DGAPA
grant IN103508.}

\vskip 1.5pc
\begin{abstract}
The starting point for this work is an identity that relates the number of
minimal matrices with prescribed 1-marginals and coefficient sequence
to a linear combination of Kronecker coefficients.
In this paper we provide a bijection that realizes combinatorially
this identity.
As a consequence we obtain an algorithm that to each minimal matrix
associates a minimal component, with respect to the dominance order,
in a Kronecker product, and a combinatorial description
of the corresponding Kronecker coefficient in terms of minimal matrices and
tableau insertion.
Our bijection follows from a generalization of the dual RSK correspondence
to 3-dimensional binary matrices, which we state and prove.
With the same tools we also obtain a generalization of the RSK correspondence
to 3-dimensional integer matrices.

\end{abstract}

\section{Introduction}\hfill

Let $\la$, $\mu$, $\nu$ be partitions of a positive integer $m$
and let $\cara\la$, $\cara\mu$, $\cara\nu$ be their corresponding
complex irreducible characters of the symmetric group $\gruposim m$.
It is a long standing problem to give a satisfactory method for
computing the multiplicity
\begin{equation*} 
\coefi := \langle \kron, \cara\nu \rangle
\end{equation*}
of $\cara\nu$ in the Kronecker product $\kron$ of $\cara\la$ and
$\cara\mu$ (here $\langle \cdot, \cdot \rangle$ denotes the inner
product of complex characters).
Via the Frobenius map, $\coefi$ is equal to the
multiplicity of the Schur function
$s_\nu$ in the internal product of Schur functions
$s_\la \ast s_\mu$, namely
\begin{equation*}
\coefi = \langle s_\la \ast s_\mu, s_\nu \rangle\, ,
\end{equation*}
where $\langle \cdot, \cdot \rangle$ denotes the scalar product of
symmetric functions.

It is not difficult to see (Lemma~\ref{lema:extremas}) that the
multiplicities of extremal (minimal or maximal) components of $\kron$
with respect to the dominance order of partitions can be described
combinatorially.
In general, the farther away a component in a Kronecker product is
from the extremal components, the harder it is to compute.
Therefore it is natural to try to determine extremal components
in a Kronecker product.

These components were investigated for the first time in~\cite{jaco},
where a connection between minimal components and discrete tomography
was discovered.
There it was shown that the existence of a minimal matrix with
row sum vector $\la$, column sum vector $\mu$ and $\pi$-sequence $\nu$
(see Section~\ref{sec:minikron} for the definitions) imply the vanishing
of ${\sf k}(\alpha, \beta, \gamma)$ for all $\alpha\mayori \la$,
$\beta\mayori\mu$ and $\gamma\menore\nu$.
It was also shown that if there exists a minimal matrix
with row sum vector $\la$, column sum vector $\mu$ and $\pi$-sequence $\nu$,
then the number $\cardimatrizcon$ of all such matrices
satisfies identity~(\ref{ecua:mini-kron}), namely,
\begin{equation*}
\cardimatrizcon
=\sum_{\alpha\mayori\la,\,\beta\mayori\mu} \kostka\alpha\lambda\kostka\beta\mu \,
\coefili \alpha\beta\nu \, .
\end{equation*}
Thus any minimal matrix on the left contributes, up to a constant, to some Kronecker
coefficient $\coefili \alpha\beta\nu$ on the right.
And in this situation, $\cara\nu$ is a minimal component of $\kronli \alpha\beta$.

In this paper we give an explicit one-to-one correspondence
(Theorem~\ref{teor:rsk-dual-tres}) that extends the dual RSK correspondence to
3-dimensional binary matrices.
By viewing each 2-dimensional integer matrix as a 3-dimensional binary matrix,
our correspondence yields a combinatorial realization of
identity~(\ref{ecua:mini-kron}).
Thus, it provides a neat way to determine the contribution of each minimal matrix
to a Kronecker coefficient  $\coefili \alpha\beta\nu$.
More precisely, we obtain a combinatorial description of the Kronecker
coefficients in~(\ref{ecua:mini-kron}) in terms of minimal matrices and
tableau insertion (Theorem~\ref{teor:mini-kron}).

In order to explain our correspondence we need the following notation.
Let $\lambda=\vector\lambda p$, $\mu=\vector\mu q$, $\nu=\vector\nu r$ be
compositions of $n$, that is, vectors with nonnegative integer coordinates
whose sum is $n$.
Let $\matriztresest$ denote the set of all 3-dimensional matrices
with entries equal to 0 or 1 and 1-marginals (plane sums) $\lambda$, $\mu$,
$\nu$ namely, binary matrices $A = (a_{ijk})$ such that their entries satisfy
\begin{equation*} 
\sum_{j,k} a_{ijk}=\lambda_i,\quad \sum_{i,k} a_{ijk}=\mu_j\quad
\textrm{and} \quad  \sum_{i,j}    a_{ijk}=\nu_k
\end{equation*}
for all $\in \conjunto p$, $j \in \conjunto q$  and $k \in \conjunto r$
(here $\conjunto n := \{1,2,\dots, n\}$).

We give a one-to-one correspondence  between $\matriztresest$
and the set of all triples $(Q, P, (T, S))$, where $Q$ and $P$ are semistandard
tableaux, not necessarily of conjugate shape, of content $\lambda$ and $\mu$
respectively, and $(T, S)$ is a pair of Littlewood-Richardson multitableaux of
conjugate content and type $\nu$, such that $T$ has the shape of $Q$ and $S$
has the shape of $P$.
This correspondence makes use of the dual RSK correspondence on each
horizontal slice of $A$.
When $\nu = ( n )$, $A$ has exactly one horizontal
slice and our correspondence is just the usual dual RSK correspondence.
In this case the pair $(T, S)$ is completely determined by $(Q,P)$.
In this sense, our bijection is a generalization
of the dual RSK correspondence to 3-dimensional $(0,1)$-matrices.
Its proof uses the dual RSK correspondence and the Littlewood-Richardson rule.

Let now $\matriztres$ denote the set of all 3-dimensional matrices
with nonnegative integer entries and 1-marginals $\lambda$, $\mu$, $\nu$.
With no extra effort we are also able to give a one-to-one correspondence between
$\matriztres$ and the set of all triples $(Q, P, (T, S))$, where $Q$ and $P$ are
semistandard tableaux, not necessarily of the same shape, of content $\lambda$ and
$\mu$, respectively, and $(T, S)$ is a pair of Littlewood- Richardson multitableaux
of the same content and type $\nu$, such that $T$ has the shape of $Q$ and $S$ has
the shape of $P$.
This is also a generalization of the RSK correspondence to 3-dimensional matrices
with nonnegative integer entries.
These bijections have already been presented without proofs in~\cite{fpsac}.

The paper is organized as follows:
In Section~\ref{sec:rsk} we explain, using character theory, how we arrive
to our main correspondences.
In Section~\ref{sec:lrm} we introduce certain pairs of sets of
Litlewood-Richardson multitableaux which are needed in our main
correspondences.
We also show that they describe combinatorially the multiplicities of
extremal components in Kronecker products (Lemma~\ref{lema:extremas}).
Section~\ref{sec:principal} contains the statement of the main theorems
(Theorems~\ref{teor:rsk-tres} and~\ref{teor:rsk-dual-tres}) and their proofs.
These theorems describe our two one-to-one correspondences that extend respectively
the RSK correspondence and its dual to 3-dimensional matrices.
Their proofs use Thomas' proof of the Littlewood-Richardson rule~\cite{thom}.
In Section~\ref{sec:ejemplo} we provide a detailed example of how the correspondence
from Theorem~\ref{teor:rsk-dual-tres} is defined.
We start Section~\ref{sec:minikron} by recalling the definition of minimal matrix
and some related results.
Then we show how Theorem~\ref{teor:rsk-dual-tres} yields a combinatorial
realization of identity~(\ref{ecua:mini-kron}).
This establishes an explicit link between minimal matrices and some Kronecker
coefficients.
As a consequence we obtain our third main result (Theorem~\ref{teor:mini-kron}):
a new combinatorial description of the multiplicities
of some minimal components in Kronecker products in terms of
minimal matrices and tableau insertion.
This is illustrated at the end of the section with an example.
Finally, Section~\ref{sec:dtk} contains, for the benefit of the reader, a brief
summary of how some notions from discrete tomography apply to Kronecker products.

\section{Matrices and RSK correspondences}\label{sec:rsk}\hfill

In this section we motivate, using character theory, our main correspondences
(presented in Section~\ref{sec:principal}).
We start with a known formula that relates the number of integral matrices
with prescribed 1-marginals to certain inner products of
characters~(\ref{ecua:clasica}).
This formula yields a second one for which the RSK-correspondence is a
combinatorial realization~(\ref{ecua:rsk}).
A similar approach can be carried out for binary matrices with
prescribed 1-marginals and the dual RSK-correspondence
(equations~(\ref{ecua:clasica-dual}) and~(\ref{ecua:rskd})).

In a similar way, but now working with 3-dimensional matrices with
prescribed 1-marginals we obtain
formulas~(\ref{ecua:rsk-tres})--(\ref{ecua:lirid}),
which suggest generalizations of the RSK and the dual RSK
correspondences to dimension~3.

Let $\permu\nu={\rm Ind}_{S_\nu}^{S_n}({\rm 1}_\nu)$ denote the character
of $S_n$ induced from the trivial character ${\rm 1}_\nu$ of the Young subgroup
$S_\nu$ associated to $\nu$.
Recall that for any $\gamma \particion n$, $\kostka\gamma\nu$
denotes the number of semistandard Young tableaux of shape $\gamma$
and content $\nu$.
Throughout this paper we will make frequent use of {\em Young's rule}, which
expresses $\permu \nu$ as a linear combination of irreducible characters,
namely
\begin{equation}
\permu\nu =\sum_{\gamma\particion n} \kostka\gamma\nu\cara\gamma .
\label{ecua:young}
\end{equation}
We use the notation $\alpha \mayori \beta$ to indicate that $\alpha$ is
bigger or equal than $\beta$ in the dominance order of partitions, and
$\alpha\mayore \beta$ if $\alpha\mayori\beta$ and $\alpha\neq\beta$.
This partial order is of interest to us because $\kostka \alpha\beta>0$ if and
only if $\alpha \mayori \beta$.

Given a matrix $A = (a_{ij})$ of size $p \times q$ we denote by
$\renglon A$ the {\bf row sum} vector of $A$ and by $\columna A$
the {\bf column sum} vector of $A$, that is, $\renglon A = \vector rp$,
where $r_i = \sum_j a_{ij}$ and $\columna M = \vector cq$, where
$c_j = \sum_i a_{ij}$.
The compositions $\renglon A$ and $\columna A$ are also called
the {\bf 1-marginals} of $A$.
Given $\la$, $\mu$ compositions of $n$, we denote by
$\matriz$ the set of all matrices $A=(a_{ij})$ with nonnegative
integer entries and 1-marginals $\la$, $\mu$,
and by $\numatriz$ its cardinality.
It is well known that $\numatriz$ can be described as an inner product
involving permutation characters and the trivial character
(see~\cite[Thm.~15]{col}, \cite[Cor.~3.1]{sna},
\cite[Thm.~1]{dfr}, \cite[6.1.9]{ker} or~\cite[7.9.1]{stan}):

\begin{equation}\label{ecua:clasica}
\numatriz = \bili{\permu\la\otimes\permu\mu}{\cara{(n)}}\, .
\end{equation}

If we expand $\permu\la$ and $\permu\mu$ as a linear combination
of irreducible characters by Young's rule~(\ref{ecua:young}),
we obtain the following identity
\begin{equation}\label{ecua:rsk}
\numatriz = \sum_{\alpha\mayori\la,\, \beta \mayori\mu} \kostka \alpha \la \kostka\beta\mu
\bili{\cara\alpha\otimes \cara\beta}{\cara{(n)}}
= \sum_{\sigma\mayori\la,\,\mu}\kostka\sigma\la \kostka\sigma\mu\, .
\end{equation}
The RSK correspondence is a combinatorial realization of this identity.

Similarly, let $\matrizd$ denote the set of all binary matrices
(matrices whose entries are either zeros or ones)
in $\matriz$ and let $\numatrizd$ denote its cardinality.
It is also known (see~\cite[Thm.~16]{col}, \cite[Cor.~7.1]{sna},
\cite[Thm.~2]{dfr} or \cite[6.1.9]{ker}) that
\begin{equation}\label{ecua:clasica-dual}
\numatrizd = \bili{\permu\la\otimes\permu\mu}{\cara{(1^n)}}\, .
\end{equation}
Applying Young's rule~(\ref{ecua:young}) we obtain the identity:
\begin{equation} \label{ecua:rskd}
\numatrizd = \sum_{\mu^{\,\prime} \mayori\sigma\mayori\la}
\kostka\sigma\la \kostka {\sigma^{\,\prime}} \mu\, .
\end{equation}
The dual RSK correspondence is a combinatorial realization of this identity.
Let us observe that identities~(\ref{ecua:clasica}) and~(\ref{ecua:clasica-dual})
can also be written in terms of inner products of symmetric functions.

There are similar results for $n$-dimensional matrices due to Snapper.
Let $\numatres$ denote the cardinality of $\matriztres$ and
$\numatresd$ denote the cardinality of $\matriztresd$.
Then Snapper showed (see Theorem~3.1 in~\cite{sna}) that
\begin{equation*} \label{ecua:snapper}
\numatres = \bili{\permu\la\otimes\permu\mu\otimes\permu\nu}{\cara{(n)}}\, .
\end{equation*}
Applying again Young's rule~(\ref{ecua:young}) we get
\begin{equation}\label{ecua:rsk-tres}
\numatres = \sum_{\alpha\mayori\la,\,\beta\mayori\mu,\,\gamma\mayori\nu}
\kostka\alpha\la \kostka\beta\mu \kostka\gamma\nu\, {\sf k}(\alpha, \beta, \gamma) \, .
\end{equation}

It is a natural question to ask for a one-to-one
correspondence between 3-dimensional matrices and some
triples of semistandard tableaux that generalizes the RSK correspondence.
Due to the nature of Kronecker coefficients, formula~(\ref{ecua:rsk-tres})
shows that this is not possible.
This formula also shows that an appropriate correspondence
(neither injective nor surjective) that associates
to a given 3-dimensional matrix a triple of semistandard tableaux would
yield a combinatorial description of Kronecker coefficients.

Similarly for binary matrices Snapper showed (see Theorem~7.1 in~\cite{sna}) that

\begin{equation*} \label{ecua:sanpperd}
\numatresd = \bili{\permu\la\otimes\permu\mu\otimes\permu\nu}{\cara{(1^n)}}\, ,
\end{equation*}
and, by Young's rule~(\ref{ecua:young}), we have
\begin{equation} \label{ecua:rskd-tres}
\numatresd = \sum_{\alpha\mayori\la,\,\beta\mayori\mu,\,\gamma\mayori\nu}
\kostka\alpha\la \kostka\beta\mu \kostka\gamma\nu\, {\sf k}(\alpha, \beta, \gamma^\prime) \, .
\end{equation}

Formulas~(\ref{ecua:rsk-tres}) and~(\ref{ecua:rskd-tres})
give us a hint of how generalizations of
the RSK correspondence and its dual to 3-dimensional matrices should be.
Let us stress that the ultimate goal is not to extend the RSK correspondence or its
dual to dimension 3 for its own sake, but to use such an extension to obtain combinatorial
descriptions of Kronecker coefficients.

Caselli~\cite[\S 4]{cas} has found some properties such generalizations of the RSK
correspondence and its dual should satisfy in order to yield combinatorial
descriptions of Kronecker coefficients.

A variation of these ideas leads to a more modest but more realistic approach:
We apply Young's rule to only two factors.
By doing so, we get formulas in which all terms have a
combinatorial description (see Lemma~\ref{lema:liri-cara}):

\begin{equation} \label{ecua:liri}
\numatres = \sum_{\alpha\mayori\la,\,\beta\mayori\mu} \kostka\alpha\la \kostka\beta\mu\,
\bili{\cara\alpha\otimes \cara\beta\otimes \permu\nu}{\cara{(n)}} \, .
\end{equation}
and
\begin{equation} \label{ecua:lirid}
\numatresd = \sum_{\alpha\mayori\la,\,\beta\mayori\mu} \kostka\alpha\la \kostka\beta\mu\,
\bili{\cara\alpha\otimes \cara\beta\otimes \permu\nu}{\cara{(1^n)}} \, .
\end{equation}

Our main results in Section~\ref{sec:principal} give combinatorial realizations of
formulas~(\ref{ecua:liri}) and~(\ref{ecua:lirid}).

\section{Littlewood-Richardson multitableaux}\label{sec:lrm}\hfill

In this section we introduce two kinds of sets of pairs of
Littlewood-Richardson multitableaux.
They will be used in the combinatorial realizations of
formulas~(\ref{ecua:liri}) and~(\ref{ecua:lirid}).
They can also be viewed as combinatorial approximations of Kronecker coefficients
(equations~(\ref{ecua:lr-kron}) and~(\ref{ecua:lrd-kron})).
In addition we deal with extremal (minimal and maximal) components $\cara\nu$ of
Kronecker products $\kronli \alpha\beta$ with respect to the dominance order
of partitions and observe that their Kronecker coefficients are
combinatorially described by those sets of pairs
(Lemma~\ref{lema:extremas}).

For the undefined terms we refer the reader to~\cite{ful,sag,stan}.

Let $\alpha$ be a partition of $n$ and
$\nu = \vector \nu r$ be a composition of $n$, then a sequence
$T = \vector Tr$ of tableaux is called a {\bf Littlewood-Richardson
multitableau} of {\bf shape} $\alpha$ and {\bf type} $\nu$ if there exists a
sequence of partitions
\[
\emptyset = \alpha(0) \subseteq \alpha(1) \subseteq \cdots
\subseteq \alpha(r) = \alpha
\]
such that $T_i$ is a Littlewood-Richardson tableau of shape
$\alpha(i)/\alpha(i -1)$ and size $\nu_i$ (the number of squares of
$T_i$ is $\nu_i$) for all $i\in \conjunto r$.
If each $T_i$ has content $\rho(i)$, then we say that $T$ has content
$(\rho(1),\dots, \rho(r))$.
Note that, since $T_i$ is a Littlewood-Richardson tableau, $\rho(i)$
is a partition of $\nu_i$.
See Section~\ref{sec:ejemplo} for an example.

Given partitions $\alpha$, $\beta$ of $n$ and $\nu$ a composition of $n$,
we denote by $\lr$ the set of all pairs $(T,S)$ of Littlewood-Richardson
multitableaux of shape $(\alpha,\beta)$ and type $\nu$ such that $S$ and
$T$ have the same content and by $\liri\alpha\beta\nu$ its cardinality.
Similarly, let $\lrd$ denote the set of all pairs $(T,S)$ of
Littlewood-Richardson multitableaux of shape $(\alpha,\beta)$,
type $\nu$ and {\bf conjugate} content, that is, if $T$ has content
$(\rho(1),\dots,\rho(r))$, then $S$ has content
$(\rho(1)^\prime,\dots, \rho(r)^\prime)$ and by $\lirid\alpha\beta\nu$
its cardinality.
Here $\rho^{\prime}$ denotes the partition conjugate to $\rho$.
We have

\begin{lema} \label{lema:liri-cara}
Let $\alpha$, $\beta$ be partitions of $n$ and let $\nu$ be a
composition of $n$.
Then

(1) $\liri \alpha\beta\nu= \bili{\cara\alpha\otimes \cara\beta\otimes \permu\nu}{\cara{(n)}}$.

(2) $\lirid \alpha\beta\nu= \bili{\cara\alpha\otimes \cara\beta\otimes \permu\nu}{\cara{(1^n)}}$.
\end{lema}

The proof of this lemma follows from Frobenius reciprocity and
the Littlewood-Richardson rule.
Part (1) of the Lemma appears implicitly in~\cite[2.9.17]{jake} and explicitly
in~\cite{estab1, estab2}.
Part (2) is similar and appears in an equivalent form in identity~(8) in~\cite{jaco}.

A component $\cara\nu$ of $\cara\alpha\otimes \cara\beta$ is called
{\bf maximal} if for all $\gamma \mayore \nu$ one has
$\coefili \alpha\beta\gamma =0$, and it is called {\bf minimal}
if for all $\gamma \menore \nu$ one has
$\coefili \alpha\beta\gamma =0$.
Minimal components were studied for the first time in~\cite{jaco}.
Since conjugation is an order-reversing involution in the set of partitions of
$n$ under the dominance order (see~\cite[6.1.18]{ker}) and since Kronecker
coefficients satisfy the symmetry $\coefili \alpha\beta {\gamma^\prime} =
\coefili \alpha{\beta^\prime} \gamma$, the study of maximal
components can be reduced to the study of minimal components.
An algorithm for computing extremal components in the lexicographic order
of partitions is given in~\cite{cm}.

It follows directly from Young's rule~(\ref{ecua:young}) that

\begin{equation} \label{ecua:lr-kron}
\liri \alpha\beta\nu = \sum_{\gamma\mayori\nu} \kostka \gamma\nu\,
\coefili\alpha\beta\gamma
\end{equation}
and
\begin{equation} \label{ecua:lrd-kron}
\lirid \alpha\beta{\nu^{\,\prime}} = \sum_{\gamma\menori\nu}
\kostka {\gamma^\prime}{\nu^{\,\prime}}\,
\coefili\alpha\beta\gamma
\end{equation}

So that, the numbers $\liri \alpha\beta\nu$ and $\lirid \alpha\beta{\nu^{\,\prime}}$
are combinatorial approximations of Kronecker coefficients.
Moreover, for extremal (maximal or minimal) components they coincide with
Kronecker coefficients:

\begin{lema} \label{lema:extremas}
Let $\cara\nu$ be a component of $\kronli\alpha \beta$.
Then

(1) $\cara\nu$ is a maximal component of $\kronli\alpha\beta$ if and only if
$\coefili \alpha\beta\nu= \liri\alpha\beta\nu$.

(2) $\cara\nu$ is a minimal component of $\kronli\alpha\beta$ if and only if
$\coefili \alpha\beta\nu = \lirid\alpha\beta{\nu^{\,\prime}}$.
\end{lema}

The proof of this lemma is straightforward.
It follows from~(\ref{ecua:lr-kron}) and~(\ref{ecua:lrd-kron}).
Part (2) is already implicit in Corollary~3.3.2 from~\cite{jaco}.

\section{Combinatorial realizations}\label{sec:principal}\hfill

In this section we give explicit bijections that are combinatorial realizations
of formulas~(\ref{ecua:liri}) and~(\ref{ecua:lirid}).

Let $\lambda$, $\mu$, $\nu$ be compositions of $n$.
For any partition $\alpha$ of $n$ let $\kostkasf \alpha\lambda$ denote
the set of all semistandard tableaux of shape $\alpha$ and content $\lambda$.
Our main results are

\begin{teor} \label{teor:rsk-tres}
There is a one-to-one correspondence between the set $\matriztres$
of 3-dimensional matrices with nonnegative integer coefficients that
have 1-marginals $\lambda$, $\mu$, $\nu$ and the set of triples
$\coprod_{\alpha\mayori\la,\,\beta\mayori\mu} \kostkasf \alpha\lambda\times\kostkasf\beta\mu
\times\lr$.
\end{teor}

\begin{teor} \label{teor:rsk-dual-tres}
There is a one-to-one correspondence between the set $\matriztresest$
of 3-dimensional binary matrices that have 1-marginals
$\lambda$, $\mu$, $\nu$ and the set of triples
$\coprod_{\alpha\mayori\la,\,\beta\mayori\mu} \kostkasf\alpha\lambda\times\kostkasf\beta\mu
\times\lrest$.
\end{teor}

The correspondences of the previous theorems will be
given as compositions of three bijections.
The first one is tautological, the second is given by a correspondence
between matrices and pairs of tableaux, such as the RSK or the dual RSK
correspondence, applied simultaneously several times, and the third
is a consequence of a bijection given by G.~P. Thomas in~\cite{thom}
for his proof of the Littlewood-Richardson rule.

Theorem~\ref{teor:rsk-tres} follows directly from the first,
the second and the third bijections given below, while
Theorem~\ref{teor:rsk-dual-tres}
follows from the first, the second and the third dual bijections.

In the statement of the bijections we use the following notation:
If $T$ is a semistandard tableau, $\forma T$ denotes its {\bf shape},
$\contenido T$ its {\bf content} and $|T|$ its {\bf size}.

We also let $p$, $q$ and $r$ denote the number of parts of
$\la$, $\mu$ and $\nu$, respectively.

\bigskip
{\em First bijection.}
There is a one-to-one correspondence between the set of matrices $\matriztres$
and the set of $r$-tuples $\vector Ar$ of matrices with nonnegative
integer coefficients of size $p \times q$ such that
\begin{equation}\label{biyec1}
\begin{gathered}
\sum_{k=1}^r \renglon{A_k} =\lambda,\qquad \sum_{k=1}^r\columna{A_k}=\mu,\\
\text{sum of the entries of } A_k=\nu_k,\  k\in \conjunto r.
\end{gathered}
\end{equation}

To construct this bijection we split $A \in\matriztres$ into its {\em level
matrices} $A^{(k)} = \left( a_{ij}^{(k)}\right)$, $k\in \conjunto r$, where
$a_{ij}^{(k)}= a_{ijk}$. Then
\[
A\asocia  \left(A^{(1)},\dots,A^{(r)}\right)
\]
is the desired bijection.

\bigskip
{\em First dual bijection.}
There is a one-to-one correspondence between the set of binary matrices
$\matriztresest$ and the set of $r$-tuples $\vector Ar$ of binary matrices
satisfying~(\ref{biyec1}).
This correspondence is the restriction of the first bijection to
$\matriztresest$.

\bigskip
{\em Second bijection.}
There is a one-to-one correspondence between the set
of $r$-tuples $\vector Ar$ of matrices with nonnegative integer coefficients
of size $p \times q$ satisfying (\ref{biyec1}) and the set of pairs
$(\vector Pr,\vector Qr)$ of $r$-tuples of semistandard tableaux such that
\begin{equation}\label{biyec2}
\begin{gathered}
\sum_{k=1}^r \contenido{Q_k} = \lambda, \qquad
\sum_{k=1}^r \contenido{P_k} = \mu, \\
\forma{P_k} = \forma{Q_k}\ {\rm and} \ \vert \forma{P_k}\vert = \nu_k,
\ k\in \conjunto r.
\end{gathered}
\end{equation}

In order to establish this bijection we choose any one-to-one
correspondence between matrices $M$ with nonnegative integer coefficients and
pairs $(P, Q)$ of semistandard tableau of the same shape such that
$\contenido P = \columna M$ and $\contenido Q = \renglon M$.
Examples of such correspondences are the RSK correspondence
\cite{knu}, \cite[4.1]{ful}, \cite[4.8]{sag} and \cite[7.11]{stan}, and the
Burge correspondence \cite[p. 198]{ful}.
Then the bijection is as follows: For any $r$-tuple $\vector Ar$ of matrices
satisfying (\ref{biyec1}), let $(P_k, Q_k)$ be the pair associated to $A_k$
under the chosen correspondence, then
\[
\vector Ar\asocia  (\vector Pr,\vector Qr)
\]
is the desired bijection.

\bigskip
{\em Second dual bijection.}
There is a one-to-one correspondence between the
set of $r$-tuples $\vector Ar$ of binary matrices of size
$p\times q$ satisfying (\ref{biyec1}) and the set of pairs
$(\vector Pr,\vector Qr)$ of $r$-tuples of semistandard tableaux such that
\begin{equation}\label{biyec2d}
\begin{gathered}
\sum_{k=1}^r \contenido{Q_k} = \lambda, \qquad
\sum_{k=1}^r \contenido{P_k} = \mu, \\
\forma{P_k} = \forma{Q_k}^\prime\ {\rm and} \ \vert \forma{P_k}\vert = \nu_k,
\ k\in \conjunto r.
\end{gathered}
\end{equation}

In order to establish this bijection we choose any one-to-one correspondence
between binary matrices $M$ and pairs $(P, Q)$ of semistandard tableaux of
conjugate shape such that $\contenido P = \columna M$ and
$\contenido Q = \renglon M$.
Examples of such correspondences are the dual RSK correspondence
\cite{knu}, \cite [p. 203]{ful}, \cite[4.8]{sag} and \cite[7.14]{stan}, and the
dual of the Burge correspondence \cite[p. 205]{ful}.
The construction of this bijection is analogous to the one of the second
bijection.

\bigskip
The two remaining bijections are based on the following result due to
G.~P. Thomas (see the Corollary in page~29 from~\cite{thom}).
There he stated it for $r=2$, but the generalization for arbitrary $r$
is straightforward.
We present his result in a slightly different form.

\begin{teor} \label{teor:thomas}
There is a one-to-one correspondence between the set of all
$r$-tuples $\vector Pr$ of semistandard tableaux and the set of pairs $(P, S)$
such that $P$ is a semistandard tableau and $S$ is a Littlewood-Richardson
multitableau of shape $\forma P$. Moreover, under this correspondence
\[
\contenido P = \sum_{k=1}^r \contenido{P_k} \quad {\rm and} \quad
\contenido S = \left( \forma{P_1}, \dots, \forma{P_r}\right).
\]
\end{teor}

\bigskip
{\em Third bijection.}
There is a one-to-one correspondence between the set of pairs of $r$-tuples
$(\vector Pr,\vector Qr)$ of semistandard tableaux satisfying (\ref{biyec2})
and the set
\begin{equation*}
\coprod_{\alpha\mayori\la,\,\beta\mayori\mu} \kostkasf\alpha\lambda\times
\kostkasf\beta\mu\times \lr.
\end{equation*}

The third bijection is as follows:
Let $(\vector Pr,\vector Qr)$ be a pair of $r$-tuples satisfying
(\ref{biyec2}), and let $(P,S)$, respectively $(Q,T)$, be the pair
corresponding to $\vector Pr$, respectively $\vector Qr$, under the
bijection of Theorem~\ref{teor:thomas}.
Then
\[
(\vector Pr,\vector Qr) \asocia (Q,P,(T,S))
\]
is the desired bijection.

\bigskip
{\em Third dual bijection.}
There is a one-to-one correspondence between the set of pairs of $r$-tuples
$(\vector Pr,\vector Qr)$ of semistandard tableaux satisfying (\ref{biyec2d})
and the set
\begin{equation*}
\coprod_{\alpha\mayori\la,\,\beta\mayori\mu} \kostkasf\alpha\lambda\times
\kostkasf\beta\mu\times \lrest.
\end{equation*}
This bijection is constructed in a similar way as the third.

\begin{obse}
{\em The correspondence in Theorem~\ref{teor:rsk-tres} satisfies, when we use the
RSK correspondence
in the second bijection of the construction, an obvious symmetry, which is
inherited from the symmetry of the RSK correspondence, namely if $A=(a_{ijk})$
of size $p\times q\times r$ corresponds to $(Q,P,(T,S))$, then its
{\em transpose} $A^t=(a_{jik})$ of size $q\times p\times r$ corresponds to
$(P,Q,(S,T))$.
Also the symmetry theorem given in~\cite[p. 205]{ful} is inherited by the
construction given in the proof of Theorem~\ref{teor:rsk-dual-tres}.}
\end{obse}

\section{An example}\label{sec:ejemplo}\hfill

In this section we explain how Thomas' bijection is defined
and give an example of the correspondence in Theorem~\ref{teor:rsk-dual-tres}

For the definition of column insertion we refer the reader to~\cite{ful, sag, stan}.
Let $x \rightarrow T$ denote the result of column inserting $x$ in a
semistandard tableau $T$.
For any partition $\gamma$, we denote by $\cano\gamma$ the unique
semistandard tableau of shape $\gamma$ and content $\gamma$.
Besides, given a semistandard tableau $T$, let $\wcol T$ denote the {\bf column  word}
of $T$, that is, the word obtained from $T$ by reading its entries from bottom to top
(in english notation), in successive columns, starting in the left column and moving to
the right.
For example $\wcol{\cano{3,2,1}} = 321211$.

Thomas' bijection is as follows:
Let $\vector Pr$ be an $r$-tuple of semistandard tableau, and let
$\gamma(k) =\forma{P_k}$, $k\in \conjunto r$.
The pair $(P,S)$ associated to $\vector Pr$ is constructed as follows.
Let $P^{(1)} =P_1$ and $S_1=\cano{\gamma(1)}$.
Then we define $P^{(k+1)}$ and $S_{k+1}$ inductively:
Let $\wcol {P_{k+1}}=v_m\cdots v_1$ and
$\wcol {\cano{\gamma(k+1)}}=u_m\cdots u_1$.
Then $P^{(k+1)}$ is obtained by column inserting  $v_1, \dots, v_m$ in
$P^{(k)}$, that is,
\begin{equation*}
P^{(k+1)}=\insertar vm{P^{(k)}},
\end{equation*}
and $S_{k+1}$ is the tableau obtained by placing  $u_1, \dots, u_m$
successively in the new boxes.
Let $P = P^{(r)}$ and $S =\vector Sr$.
Then $P$ is a semistandard tableau, $S$ is a Littlewood-Richardson
multitableau, $\forma P = \forma S$, $\contenido
S=(\gamma(1),\dots,\gamma(r))$, and $\contenido P=\sum_{k=1}^r
\contenido{P_k}$.
Note that $P=P_r \boldsymbol{\cdot} \cdots \boldsymbol{\cdot} P_1$, the product of
tableaux as defined in Fulton's book~\cite{ful}.

\bigskip
We conclude this section with an illustration of the bijection described
in Theorem~\ref{teor:rsk-dual-tres}; for this we use the dual RSK
correspondence in the second bijection of the construction.
Let $A$ be the following 3-dimensional matrix of zeroes and ones of size
$4 \times 5 \times 3$.
\begin{equation*}
\left[
\begin{matrix}
0 & 1 & 1 & 1 & 0 \\
1 & 1 & 0 & 0 & 1 \\
0 & 0 & 1 & 1 & 0 \\
1 & 0 & 1 & 0 & 0
\end{matrix}
\right]
\qquad
\left[
\begin{matrix}
1 & 0 & 1 & 0 & 1 \\
1 & 0 & 0 & 0 & 0 \\
1 & 1 & 0 & 1 & 0 \\
0 & 1 & 0 & 0 & 0
\end{matrix}
\right]
\qquad
\left[
\begin{matrix}
0 & 1 & 0 & 1 & 1 \\
1 & 1 & 1 & 0 & 0 \\
0 & 0 & 0 & 0 & 0 \\
1 & 0 & 0 & 0 & 0
\end{matrix}
\right]
\end{equation*}
It has 1-marginals $\lambda = (9,7,5,4)$, $\mu = (7,6,5,4,3)$ and
$\nu = (10,8,7)$.
The triple $(Q, P, (T, S))$ corresponding to $A$
is constructed as follows:
To each of the three level matrices corresponds, under the dual RSK
correspondence, a pair of semistandard tableaux of conjugate shape
\[
(P_1,Q_1)=
\left(\
\raisebox{1.7em}[2em][3em]
{\tableau{
1 & 1 & 2 \\
2 & 3 & 3 \\
3 & 4 \\
4 & 5 \\}}
\ ,\
\raisebox{1.7em}
{\tableau{
1 & 1 & 1 & 2 \\
2 & 2 & 3 & 3 \\
4 & 4 \\}}\
\right)
\]
\begin{equation*}
(P_2,Q_2)=
\left(\
\raisebox{.95em}[2.8em]
{\tableau{
1 & 1 & 1 \\
2 & 2 & 3 \\
4 & 5 \\}}
\ ,\
\raisebox{.95em}
{\tableau{
1 & 1 & 1 \\
2 & 3 & 3 \\
3 & 4 \\}}\
\right)
\qquad
(P_3,Q_3)=
\left(\
\raisebox{.95em}[2.8em]
{\tableau{
1 & 1 & 2 \\
2 & 4 \\
3 & 5 \\}}
\ ,\
\raisebox{.95em}
{\tableau{
1 & 1 & 1 \\
2 & 2 & 2 \\
4 \\}}\
\right)
\end{equation*}
Then $(Q,T)$ and $(P,S)$ are the pairs associated to $(Q_1, Q_2, Q_3)$
and $(P_1, P_2, P_3)$, respectively, under the correspondence given in
the proof of Theorem~\ref{teor:thomas}.
The pair $(Q,P)$ of semistandard tableaux is
\[
\left(\
\raisebox{2em}
{\tableau{
1 & 1 & 1 & 1 & 1 & 1 & 1 & 1 & 1 & 2 \\
2 & 2 & 2 & 2 & 2 & 2 & 3 & 3 \\
3 & 3 & 3 & 4 & 4 \\
4 & 4 \\}}
\quad ,\quad
\raisebox{2em}
{\tableau{
1 & 1 & 1 & 1 & 1 & 1 & 1 & 2 & 3 \\
2 & 2 & 2 & 2 & 2 & 3 & 3 \\
3 & 3 & 4 & 4 & 4 \\
4 & 5 & 5 \\
5 \\}}\
\right)
\]
and the pair $(T,S)$ of Littlewood-Richardson multitableaux is
\[
\left(\
\raisebox{2em}
{\tableau{
{\bf 1} & {\bf 1} & {\bf 1} & {\bf 1} & {\it 1} & {\it 1} & {\it 1} & {\tt 1}
& {\tt 1} & {\tt 1} \\
{\bf 2} & {\bf 2} & {\bf 2} & {\bf 2} & {\it 2} & {\tt 2} & {\tt 2}& {\tt 2} \\
{\bf 3} & {\bf 3} & {\it 2}  & {\it   2}  & {\it 3} \\
{\it 3} & {\tt 3} \\}}
\quad ,\quad
\raisebox{2em}
{\tableau{
{\bf 1} & {\bf 1} & {\bf 1} & {\it 1} & {\it 1} & {\it 1} & {\tt 1} & {\tt 1} & {\tt 1} \\
{\bf 2} & {\bf 2} & {\bf 2} & {\it 2} & {\it 2} & {\it 2} & {\tt 2} \\
{\bf 3} & {\bf 3} & {\it 3} & {\tt 2} & {\tt 3} \\
{\bf 4} & {\bf 4} & {\tt 3} \\
{\it 3} \\}}\
\right)
\]
The multitableau $T$ has three parts $(T_1,T_2,T_3)$; we indicated the numbers
in $T_1$ with {\bf boldface} numerals, the numbers in $T_2$ with \textit{italic}
numerals and the numbers in $T_3$ with {\tt typewriter} numerals.
Similarly for $S$.

\section{Minimal matrices and Kronecker products}\label{sec:minikron}\hfill

Minimal matrices were introduced in \cite{tova}
to characterize 3-dimensional binary matrices that are
uniquely determined by its 1-marginals.
They were used in \cite{jaco} as a tool to produce minimal
components in Kronecker products.
In this section we go an step further towards an understanding
of the relation between minimal matrices and Kronecker products.
As an application of Theorem~\ref{teor:rsk-dual-tres} we give an
algorithm that, out of a list of minimal matrices, computes
several Kronecker coefficients.
We start by recalling some definitions and results.

For any matrix $A$ with nonnegative integer entries, we denote by
$\parti A$ the weakly decreasing sequence of its entries and call it a
{\bf $\boldsymbol{\pi}$-sequence}.
Let $\la$, $\mu$, $\nu$ be three partitions of some integer $n$.
We denote by $\matrizcon$ the subset of $\matriz$ formed by all matrices
$A$ with $\parti A=\nu$, and by $\cardimatrizcon$ its cardinality.
A matrix $A$ in $\matriz$ is called {\bf minimal} (see~ \cite{tova})
if there is no other matrix $B$ in $\matriz$ such that
$\parti B\menore \parti A$.
Note that if $A\in\matrizcon$ is minimal, then all matrices in
$\matrizcon$ are minimal.
We say that $\nu$ is {\bf minimal for}{\mathversion{bold} $(\la,\mu)$}
if there is a minimal matrix in $\matrizcon$.

\begin{ejem}
{\em Let
\[
A= \left[
\begin{matrix}
0 & 3 \\
3 & 0
\end{matrix}
\right] ,\quad
B= \left[
\begin{matrix}
1 & 2 \\
2 & 1
\end{matrix}
\right],\quad
C= \left[
\begin{matrix}
2 & 1 \\
1 & 2
\end{matrix}
\right] ,\quad
D=\left[
\begin{matrix}
3 & 0 \\
0 & 3
\end{matrix} \right],
\]
then $A$, $B$, $C$ and $D$ have the same 1-marginals $\la=\mu =(3^2)$,
$\parti A = \parti D=(3^2)$ and $\parti B=\parti C= (2^2,1^2)$.
The set $\matriz$ is equal to $\{ A, B,C,D\}$,
thus $B$ and $C$ are minimal, $A$ and $D$ are not and $(2^2,1^2)$ is minimal
for $(\la,\mu)$.}
\end{ejem}

The following theorem establishes a connection between minimal matrices
and multiplicities of minimal components in Kronecker products.

\begin{teor}
If $\nu$ is minimal for $(\la,\mu)$.
Then

(i) $\coefili \alpha\beta\gamma =0$ for all $\alpha\mayori\lambda$,
$\beta\mayori\mu$, $\gamma\menore\nu$.

(ii) $\coefili \alpha\beta\nu = \lirid \alpha\beta{\nu^{\,\prime}}$ for all
$\alpha\mayori\lambda$, $\beta\mayori\mu$.

In particular, for any pair of partitions $(\alpha,\beta)$ such that
$\alpha\mayori\lambda$ and $\beta\mayori\mu$ we have that
$\cara\nu$ is a minimal component of $\kronli \alpha\beta$ if and only if
$\lirid \alpha \beta {\nu^{\,\prime}}$ is positive.
\end{teor}

Note that (i) is Proposition 3.2 in \cite{jaco},
and (ii) follows from~(\ref{ecua:lrd-kron}) and (i).
The last remark follows from Lemma~\ref{lema:extremas}.

To each matrix $A=(a_{ij})$ in $\matrizcon$ we associate a
3-dimensional matrix $\overline A=(a_{ijk})$ by
\[
a_{ijk}=
\begin{cases}
1 & \text{if $a_{ij}\le k$}, \\
0 & \text{otherwise}.
\end{cases}
\]
The correspondence $A\mapsto \overline A$ defines an injective map
\begin{equation}
G_{\lambda,\mu,\nu} : \matrizcon \flecha {\sf M}^*(\lambda,\mu,\nu^{\,\prime}).
\label{ecua:inyeccion}
\end{equation}
We have the following characterization of minimality.

\begin{prop}
\cite[Thm. 13]{dito} \label{prop:carac-mini}
Let $\la$, $\mu$, $\nu$ be partitions of $n$.
Then $\nu$ is minimal for $(\la,\mu)$ if and only if $G_{\lambda,\mu,\nu}$ is bijective.
\end{prop}

The proof given in~\cite{dito} is combinatorial.
A different proof follows immediately from Proposition~3.1 in~\cite{jaco}
and the fact that the map $G_{\lambda,\mu,\nu}$ is injective.
Then, from this proposition, identity~(\ref{ecua:lirid}) and lemmas~\ref{lema:liri-cara}
and \ref{lema:extremas}, we obtain

\begin{coro}
\cite[Cor. 3.3.2]{jaco}
Let $\nu$ be a minimal for $(\la,\mu)$.
Then
\begin{equation} \label{ecua:mini-kron}
\cardimatrizcon
=\sum_{\alpha\mayori\la,\,\beta\mayori\mu} \kostka\alpha\lambda\kostka\beta\mu \,
\lirid \alpha\beta{\nu^{\,\prime}}
=\sum_{\alpha\mayori\la,\,\beta\mayori\mu} \kostka\alpha\lambda\kostka\beta\mu \,
\coefili \alpha\beta\nu \, .
\end{equation}
\end{coro}

Let
\begin{equation*}
\Phi^*\colon \matriztresdc \flecha
\coprod_{\alpha\mayori\la,\,\beta\mayori\mu} \kostkasf\alpha\lambda\times\kostkasf\beta\mu
\times \sf{LR}^*(\alpha,\beta;\nu^{\,\prime})\, .
\end{equation*}
denote the bijection that we get from Theorem~\ref{teor:rsk-dual-tres},
when we apply in each level the dual RSK-correspondence.
Then the composition
\begin{equation} \label{ecua:matriz-tablas}
\Phi^* \circ G_{\la,\mu,\nu} \colon \matrizcon \flecha
\coprod_{\alpha\mayori\la,\,\beta\mayori\mu} \kostkasf\alpha\lambda\times\kostkasf\beta\mu
\times \sf{LR}^*(\alpha,\beta;\nu^{\,\prime})\, .
\end{equation}
is injective.

\begin{obse}
{\em Let ${\sf P}_\nu(\la,\mu)$ denote the set of plane partitions in $\matrizcon$.
There is an injective map
\begin{equation} \label{ecua:partipla-lr}
\partipla \flecha \sf{LR}^*(\alpha,\beta;\nu^{\,\prime})\, ,
\end{equation}
which was defined in the proof of Theorem~3.4 in~\cite{jaco}.
It is straightforward to verify that~(\ref{ecua:matriz-tablas}) is an
extension of~(\ref{ecua:partipla-lr}).
The starting point of this paper was the attempt to find an extension of
this injection to $\matrizcon$.
We finally managed to extend it, not only to $\matrizcon$,
but to ${\sf M}^*(\lambda,\mu,\nu^{\,\prime})$, thus getting
Theorem~\ref{teor:rsk-dual-tres}.
This was of interest because, by Lemma~\ref{lema:extremas},
$\sf{LR}^*(\alpha,\beta;\nu^{\,\prime})$ was, in some cases, a combinatorial realization of
$\coefili \alpha\beta\nu$.
Let us note that Manivel~\cite[Prop.~3.1]{mani} showed that $|\partipla|\le \coefi$
for all $\la$, $\mu$, $\nu$, thus getting a better approximation
than~(\ref{ecua:partipla-lr}).
However, this is still, in general, a weak approximation of $\coefi$.
This approximation is better when $\nu$ is minimal for $(\la,\mu)$,
see Remark~4.4 in~\cite{jaco}.}
\end{obse}

Note that if $\nu$ is minimal for $(\la,\mu)$,
(\ref{ecua:matriz-tablas}) is a bijection and provides a combinatorial
realization of~(\ref{ecua:mini-kron}).
Thus, if $\nu$ is minimal for $(\la,\mu)$ and we have the list
of all elements in $\matrizcon$ we can compute several
Kronecker coefficients.
In fact, we get a combinatorial description of these coefficients:
Let $f$, $g$, $h$ denote the components of $\Phi^* \circ G_{\la,\mu,\nu}$,
that is, for any $A\in\matrizcon$ we have
$\Phi^* \circ G_{\la,\mu,\nu}(A) = (f(A), g(A) , h(A))$.
Then

\begin{teor} \label{teor:mini-kron}
Suppose $\nu$ is minimal for $(\la,\mu)$.
Let $P$ be a semistandard tableau of shape $\alpha$ and content $\la$,
and $Q$ be a semistandard tableau of shape $\beta$ and content $\mu$.
Then
\[
\coefili \alpha\beta\nu=
\# \{ A\in\matrizcon \mid f(A)= P \text{ and } g(A)= Q\}\, .
\]
Moreover, if $\coefili \alpha\beta\nu >0$, then $\cara\nu$ is
a minimal component of $\cara\alpha\otimes \cara\beta$.
\end{teor}

We conclude this section with an illustration of Theorem~\ref{teor:mini-kron}.

\begin{ejem}
{\em Let $\la=(6,6)$, $\mu= (3,3,3,3)$.
Then, there are six minimal matrices in $\matriz$ (see Theorem~1.1 and
Lemma~4.1 in~\cite{minidosq}), namely
\begin{equation*}
A=\left[
\begin{matrix}
2 & 2 & 1 & 1  \\
1 & 1 & 2 & 2
\end{matrix}
\right]\, ,
\quad
B=\left[
\begin{matrix}
1 & 2 & 2 & 1  \\
2 & 1 & 1 & 2  \\
\end{matrix}
\right]\, ,
\quad
C=\left[
\begin{matrix}
2 & 1 & 2 & 1  \\
1 & 2 & 1 & 2  \\
\end{matrix}
\right]\, ,
\end{equation*}
\begin{equation*}
D=\left[
\begin{matrix}
2 & 1 & 1 & 2  \\
1 & 2 & 2 & 1  \\
\end{matrix}
\right]\, ,
\quad
E=\left[
\begin{matrix}
1 & 2 & 1 & 2  \\
2 & 1 & 2 & 1  \\
\end{matrix}
\right]\, ,
\quad
F=\left[
\begin{matrix}
1 & 1 & 2 & 2  \\
2 & 2 & 1 & 1  \\
\end{matrix}
\right]\, .
\end{equation*}
Let $\nu= (2^4,1^4)$ be the common $\pi$-sequence of the six matrices.
After computing $\Phi^* \circ G_{\la,\mu,\nu}$ for each matrix we get
\begin{equation*}
f(A) = \raisebox{.35em}
{\tableau{
1 & 1 & 1 & 1 & 1 & 1 & 2 & 2 \\
2 & 2 & 2 & 2 }}
\qquad \text{and} \qquad
g(A)= \raisebox{1.75em}[3em][3em]
{\tableau{
1 & 1 & 1 \\
2 & 2 & 2 \\
3 & 3 & 3 \\
4 & 4 & 4}}\ .
\end{equation*}
\begin{equation*}
f(B) = \raisebox{.35em}
{\tableau{
1 & 1 & 1 & 1 & 1 & 1 & 2\\
2 & 2 & 2 & 2 & 2}}
\qquad \text{and} \qquad
g(B)= \raisebox{1.75em}[3em][3em]
{\tableau{
1 & 1 & 1 & 2\\
2 & 2 & 3 \\
3 & 3 & 4 \\
4 & 4 }}\ .
\end{equation*}
\begin{equation*}
f(C) = \raisebox{.35em}
{\tableau{
1 & 1 & 1 & 1 & 1 & 1 & 2\\
2 & 2 & 2 & 2 & 2 }}
\qquad \text{and} \qquad
g(C) = \raisebox{1.75em}[3em][3em]
{\tableau{
1 & 1 & 1 & 3\\
2 & 2 & 2 \\
3 & 3 & 4 \\
4 & 4 }} \ .
\end{equation*}
\begin{equation*}
f(D) = \raisebox{.35em}
{\tableau{
1 & 1 & 1 & 1 & 1 & 1 & 2\\
2 & 2 & 2 & 2 & 2 }}
\qquad \text{and} \qquad
g(D) = \raisebox{1.75em}[3em][3em]
{\tableau{
1 & 1 & 1 & 4\\
2 & 2 & 2 \\
3 & 3 & 3 \\
4 & 4 }} \ .
\end{equation*}
\begin{equation*}
f(E) = \raisebox{.35em}
{\tableau{
1 & 1 & 1 & 1 & 1 & 1 \\
2 & 2 & 2 & 2 & 2 & 2 }}
\qquad \text{and} \qquad
g(E) = \raisebox{1.75em}[3em][3em]
{ \tableau{
1 & 1 & 1 & 2\\
2 & 2 & 3 & 4 \\
3 & 3 \\
4 & 4 }} \ .
\end{equation*}
\begin{equation*}
f(F) = \raisebox{.35em}
{\tableau{
1 & 1 & 1 & 1 & 1 & 1 \\
2 & 2 & 2 & 2 & 2 & 2 }}
\qquad \text{and} \qquad
g(F) = \raisebox{1.75em}[3em][3em]
{\tableau{
1 & 1 & 1 & 3 \\
2 & 2 & 2 & 4 \\
3 & 3  \\
4 & 4 }} \ .
\end{equation*}
Let $\alpha = \forma{f(B)}=(7,5)$,
$\beta = \forma{g(B)}= (4,3,3,2)$,
$\gamma = \forma{f(A)}=(8,4)$,
$\delta = \forma{g(E)}=(4,4,2,2)$.
Thus we obtain from Theorem~\ref{teor:mini-kron} that
${\sf k}(\gamma, \mu, \nu) = 1$,
${\sf k}(\alpha, \beta, \nu) = 1$ and
${\sf k}(\la, \delta, \nu) = 1$ and that $\cara\nu$ is a minimal component of
$\kronli \gamma\mu$, $\kronli\alpha\beta$ and $\kronli\la\delta$.
The remaining Kronecker coefficients in (\ref{ecua:mini-kron}) are all zero.
For example $\coefi =0$.}
\end{ejem}

\section{Discrete Tomography and Kronecker products}\label{sec:dtk}\hfill

In this brief section we show, for the benefit of the interested
reader, how some notions from discrete tomography apply to
Kronecker products.

Let $\la$, $\mu$, $\nu$ be partitions of some integer $n$.
A matrix $X\in \matriztresd$ is called a {\bf matrix of uniqueness}
if $\numatresd=1$, that is, if $X$ is the only binary matrix with
1-marginals $\la$, $\mu$, $\nu$.
This notion appears in discrete tomography
(see~\cite{bdg, flrs, fishe, grvr, dito}) and
is of interest to Kronecker products because, if $\numatresd =1$
all Kronecker coefficients but one vanish in equation~(\ref{ecua:rskd-tres}).

There is a combinatorial characterization of uniqueness
that is useful in some instances.
We need a definition in order to explain it:
A matrix $A\in\matrizcon$ is called {\bf $\boldsymbol{\pi}$-unique}
if it is the only matrix in $\matriz$ with $\pi$-sequence $\nu$, that
is, if $\cardimatrizcon=1$.
It was shown in~\cite{tova} (see also Theorem~11 in~\cite{dito})
that a matrix $X\in {\sf M}^*(\la,\mu,\nu^{\,\prime})$
is a matrix of uniqueness if and only if there is a matrix $A\in\matrizcon$
that is minimal and $\pi$-unique and such that $G_{\la,\mu,\nu}(A)=X$
(see~(\ref{ecua:inyeccion}) for the definition of $G_{\la,\mu,\nu}$).
Therefore, the existence of a matrix $A\in\matrizcon$ that is minimal and
$\pi$-unique implies that ${\sf m}^*(\la,\mu,\nu^{\,\prime})=1$.
For example, let $A$ be the $p\times q$ matrix such that all its entries are
equal to $r$.
Then $A\in\matrizcon$ where $\la =((qr)^p)$, $\mu= ((pr)^q)$ and $\nu=(r^{pq})$.
It is very easy to see that $A$ is minimal and $\pi$-unique.
Another family of matrices that are minimal and $\pi$-unique appears
in~\cite[p.~446]{tova}.
In this case $\la$, $\mu$ and $\nu$ are hooks.

Corollary~4.2 in~\cite{jaco} summarizes all consequences of uniqueness to
Kronecker coefficients.
It can be reformulated in the following way:

\begin{teor}
Let $\la$, $\mu$, $\nu$ be partitions of $n$ and let $A\in\matrizcon$.
If $A$ is minimal and $\pi$-unique, then $\cara\nu$ is a minimal component
of $\kron$, $\coefi =1$ and $\coefili \alpha\beta\gamma =0$ for all other triples
$(\alpha,\beta,\gamma)$ such that $\alpha\mayori\la$, $\beta\mayori\mu$ and
$\gamma\menori\nu$.
\end{teor}

There is still another useful tool to determine uniqueness of a matrix.
A notion of additivity for 3-dimensional binary matrices was introduced
in~\cite{flrs} and was shown to be a sufficient condition for a matrix in
$\matriztresd$ to be a matrix of uniqueness.
This notion was later translated to a version of additivity for integer
matrices (Theorem~1 in~\cite{aditivo}):
a matrix $A=(a_{ij})$ of size $p\times q$ with nonnegative integer entries
is called {\bf additive} if there exist real numbers $x_1, \dots, x_p$
and $y_1, \dots, y_q$ such that the condition
\begin{equation*}
a_{ij} > a_{kl} \Longrightarrow x_i +y_j > x_k + y_l
\end{equation*}
holds for all $i$, $j$, $k$, $l$.
Later, the obvious extension of additivity from integer to real
matrices was studied from a geometric point of view in~\cite{onva}.
Additivity for binary matrices seems to have been motivated by
a related notion for binary relations (see~\cite{fish} and the references therein).

The next result appears as Theorem~6.1 in~\cite{aditivo} and
Corollary~6.2 in~\cite{elenot}.
A geometric proof can be found in~\cite{onva}.

\begin{teor}
Any additive matrix with nonnegative integer entries is minimal and $\pi$-unique.
\end{teor}

In particular each additive matrix in $\matrizcon$ yields a minimal
component $\cara\nu$ of $\kron$ with multiplicity 1.

\begin{obse}
{\em Minimal matrices of size $2\times q$ were classified in~\cite{minidosq}.
Any plane partition of size $2 \times q$ is additive (see Proposition~4.1
in~\cite{aditivo} and Lemma~8 in~\cite{dito} for a shorter proof).
There is no general known result for minimal or additive matrices of size $3\times q$.
However, a complete set of obstructions for additivity --the so called arrow diagrams--
was given in~\cite{sava}.
There, it was also shown (Theorem~4.1) that no finite subset of such obstructions
is enough to determine additivity of an arbitrary integer matrix.

If $a_1, a_2, \dots, a_p$ and $b_1, b_2, \dots, b_q$ are nonnegative
integers then, by the very definition of additivity, the matrix
$A=(a_i + b_j)_{i\in \conjunto p,\, j\in \conjunto q}$ is additive.
Other examples of minimal or additive matrices appear in~\cite{jaco, aditivo, dito}.}
\end{obse}

\vskip 2pc
{\bf Acknowledgments}

\bigskip
Part of this work was done while the second named author stayed at the
Isaac Newton Institute during the programm {\em Symmetric Functions and
Macdonald Polynomials}.
He would like to thank the Institute for partial support and the organizers
Phil Hanlon, Ian Macdonald and Allun Morris for the very pleasant working
atmosphere.

\end{document}